# On the Estimation of Directional Returns to Scale via DEA models

Guoliang Yang[1], Wenbin Liu[2], Wanfang Shen[3] and Xiaoxuan Li[1]

[1]Institute of Policy and Management, Chinese Academy of Sciences, Beijing 100190, China

(Email: glyang@casipm.ac.cn)

[2]Kent Business School, Canterbury CT2 7PE, UK

[3]School of Mathematic and Quantitative Economics, Shandong University of Finance and Economics, Jinan 240014, China

**Abstract.** Data envelopment analysis (DEA) is one of the most commonly used methods to estimate the returns to scale (RTS) of the public sector (e.g., research institutions). Existing studies are all based on the traditional definition of RTS in economics and assume that multiple inputs and outputs change in the same proportion, which is the starting point to determine the qualitative and quantitative features of RTS of decision making units (DMUs). However, for more complex products, such as the scientific research in institutes, changes of various types of inputs or outputs are often not in proportion. Therefore, the existing definition of RTS in the framework of DEA method may not meet the need to estimate the RTS of research institutions with multiple inputs and outputs. This paper proposes a definition of directional RTS in the DEA framework and estimates the directional RTS of research institutions using DEA models. Further in-depth analysis is conducted for an illustrative example of 16 basic research institutes in Chinese Academy of Sciences (CAS) in 2010.

**1 Introduction**

In *The Quarterly Journal of Economic*s, Panzar and Willig (1977) proposed a method to determine the returns to scale (RTS) of decision-making units (DMUs) based on the production function. The estimation of RTS of DMUs using the data envelopment analysis (DEA) method was investigated first by Banker (1984) and Banker et al. (1984). Banker (1984) introduced the definition of the RTS from classical economics into the framework of the DEA method, and he used the CCR-DEA model with radial measure to estimate the RTS of evaluated DMUs. Soon after that, Banker et al. (1984) proposed the BCC-DEA model under the assumption of variable RTS and investigated how to apply the BCC-DEA model to estimate the RTS of DMUs. Thus far, in addition to the cost-based measurement of RTS of DMUs (e.g., Färe and Grosskopf, 1985; Färe et al., 1994; Seitz, 1970; Sueyoshi, 1999), DEA-based studies of DMUs' RTS can be roughly divided into four categories: (1) RTS measurement using CCR-DEA models, (2) RTS measurement using BCC-DEA models, (3) RTS measurement using FGL-DEA model and quantitative measurement of scale elasticity (SE) and (4) RTS measurement using non-radial DEA models. Please see the literature review in Section 2 for details.

The existing studies on the RTS measurement in DEA models are all based on the definition of RTS in the DEA framework made by Banker (1984), which extended the application area of DEA from relative efficiency evaluation to RTS measurement. The RTS is a classic economic concept describing the relationship between changes in the scale of production and output. The traditional definition of RTS in economics is based on the idea of measuring radial changes in outputs caused by all inputs. For example, the SE (right-hand) of 1.5 tells us that the increase of all inputs by, say, 1% corresponds to the increase of the outputs by 1.5%. Following this concept, Banker (1984) defined RTS in the DEA framework using the radial changes in outputs caused by all inputs.

In traditional industrial production, the proportion of labour and capital inputs are often fixed, so using a radial idea to define RTS in economics is practical. However, in research organisations, it often can be observed that production factors are not necessarily tied together proportionally because of the complexity of research activities and inputs change non-proportionally, as illustrated by the following Example 1.

**Example 1:** In the period of 1998-2007, both the S&T inputs and outputs increased significantly in the Chinese Academy of Sciences (CAS). We selected staff and funding as input indicators and international papers as one of the outputs. The full-time equivalent (FTE) of R&D personnel at the CAS grew from 30,611 in 1998 to 44,307 in 2007. The total funding for the CAS grew from 4,935.98 million RMB in 1998 to 17,039.71 million RMB in 2007. The number of international papers grew from 5,478 in 1998 to 24,045 in 2007. The proportions of annual changes of the two input indicators are very different as shown in Table 1.

Table 1: The changes of some indicators for the CAS from 1998 to 2007

| Year | Input indicators | Output indicators |
|---|---|---|

|  | Total funding | | R&D personnel | | International papers | |
| --- | --- | --- | --- | --- | --- | --- |
|  | Total amount (million RMB) | Growth proportion | FTE | Growth proportion | Number | Growth proportion |
| 1998 | 4,935.98 | 10.46% | 30,611 | -7.11% | 5,478 | 50.58% |
| 1999 | 5,452.06 | 30.92% | 28,436 | -1.24% | 8,249 | 11.36% |
| 2000 | 7,137.98 | 12.93% | 28,084 | -10.27% | 9,186 | 10.66% |
| 2001 | 8,060.83 | 24.98% | 25,199 | 9.71% | 10,165 | 15.49% |
| 2002 | 10,074.21 | -2.94% | 27,646 | 11.90% | 11,740 | 23.65% |
| 2003 | 9,778.10 | 24.94% | 30,937 | 12.80% | 14,516 | 8.42% |
| 2004 | 12,216.49 | 4.38% | 34,898 | 6.73% | 15,738 | 41.42% |
| 2005 | 12,751.83 | 14.12% | 37,246 | 4.47% | 22,257 | 5.98% |
| 2006 | 14,552.50 | 17.09% | 38,911 | 13.87% | 23,589 | 1.93% |
| 2007 | 17,039.71 | N/A | 44,307 | N/A | 24,045 | N/A |

Data source: Statistical Yearbook of Chinese Academy of Sciences, 1999-2008.

From Table 1, we can see that the inputs do not change proportionally. In fact, in radial measurement, the inputs did not even increase for some years. However, it is also clear that the outputs have greatly increased during the period. The traditional definition of RTS based on radial measure strays considerably from the reality of the input changes. Therefore, we need to consider directional returns to scale (directional RTS) with non-proportional changes in inputs (or outputs), as seen in Figs 1-2 below. In such a case of non-proportional change, for example, the directional right-hand SE (See definitions on directional SE and directional RTS in Section 3) of 1.5 tells us that the increase of all inputs by 1% in certain input direction corresponds to the increase of the outputs by 1.5% in certain output direction. It should be noted that the directional RTS is still based on the Pareto preference.

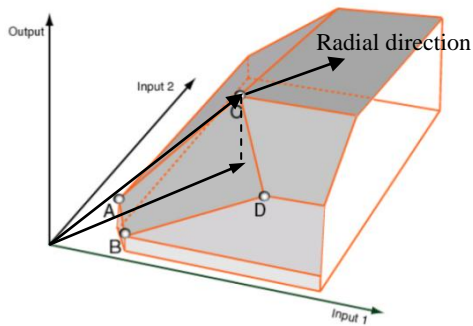 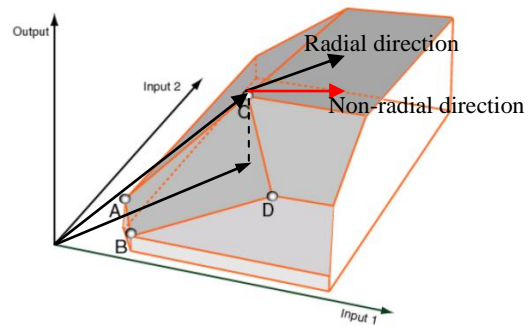

Figure 1: Traditional RTS　　　　　Figure 2: Directional RTS

This paper aims to investigate the RTS measurement of DMUs on the efficient frontier of production possibility set (PPS) in the case of inputs and outputs changing in unequal proportions, which is essentially different with the RTS measurement using non-radial models.

The rest of this paper is organised as follows. Section 2 provides the classical RTS in the DEA framework. Section 3 proposes the definitions of directional RTS and directional SE. Two approaches, the finite difference method (FDM) and upper and lower bounds method (ULBM), are proposed in Section 4 to determine the directional RTS of efficient DMUs on the efficient frontier. Section 5 provides an illustrative example of the in-depth analysis of 16 basic CAS research institutes in 2010. The last section offers the conclusion.

## 2 Classical RTS in DEA framework

We consider a set of n observations of actual production possibilities $(X_j, Y_j), j=1,...,n$. The output vector $Y_j$ can be produced from the input vector $X_j$. First, we provide the following definitions:

**Definition 1:** The PPS under the assumption of variable RTS is defined as follows.

$$PPS = \left\{(X,Y) \Big| X \geq \sum_{j=1}^{n} \lambda_j X_j, Y \leq \sum_{j=1}^{n} \lambda_j Y_j, \sum_{j=1}^{n} \lambda_j = 1, \lambda_j \geq 0, j=1,...,n \right\} \quad (1)$$

**Definition 2:** The weakly and strongly efficient frontiers of PPS can be defined as follows.

(1) Weakly efficient frontier:

$$EF_{weak} = \left\{(X,Y) \in PPS \Big| \text{there is no } (\bar{X},\bar{Y}) \in PPS \text{ such that } (-\bar{X},\bar{Y}) > (-X,Y) \right\} \quad (2)$$

(2) Strongly efficient frontier:

$$EF_{strong} = \left\{(X,Y) \in PPS \Big| \text{there is no } (\bar{X},\bar{Y}) \in PPS \text{ such that } (-\bar{X},\bar{Y}) \geq (-X,Y) \text{ and } (\bar{X},\bar{Y}) \neq (X,Y) \right\} \quad (3)$$

**Definition 3:** A supporting hyperplane of PPS can be defined as follows. If a hyperplane

$$H(V,U,\mu_0) = \{(X,Y) | U^T Y - V^T X + \mu_0 = 0\} \quad (4)$$

satisfies (1) $(X_0, Y_0) \in \{(X,Y) | U^T Y - V^T X + \mu_0 = 0\} \cap PPS$;

(2) $U^T Y - V^T X + \mu_0 \leq 0$ for all $(X,Y) \in PPS$;

(3) $(V,U) \neq (\mathbf{0},\mathbf{0})$

where $U = (u_1, u_2, ..., u_s)^T$ and $V = (v_1, v_2, ..., v_m)^T$ are vectors of multipliers, $\mu_0$ is a parameter free of sign.

Then, we say that $H(V,U,\mu_0)$ is a supporting hyperplane of PPS on the point $(X_0,Y_0)$, which is referred as $H(V,U,\mu_0)|_{(X_0,Y_0)}$.

**Definition 4:** A subset of PPS is referred to as a "Face" on the point $(X_0,Y_0)$ if there exists a supporting hyperplane $H(V,U,\mu_0)|_{(X_0,Y_0)}$ such that the subset is identical to an intersection between PPS and $H(V,U,\mu_0)|_{(X_0,Y_0)}$.

Banker (1984) defined the classical RTS based on PPS in DEA framework. As mentioned in the first section, in addition to the cost-based measurement of RTS of DMUs, DEA-based studies of DMUs' RTS can be roughly divided into four categories:

(1) RTS measurement using CCR-DEA models. The first and well-known approach to determine the RTS of DMUs is to calculate the value of $\sum_{i=1}^{n}\lambda_j^*$ in CCR-DEA models, where $\lambda_j^*$ denotes the weight of $DMU_j$ (Note: In some cases $\lambda_j^*$ may be not unique). Research efforts related to this approach can be found in Banker (1984), Chang and Guh (1991), Banker and Thrall (1992), Zhu and Shen (1995), Banker et al. (1996a, 1996b), Seiford and Zhu (1998, 1999), among others. In this approach, the type of RTS can be determined as increasing RTS, constant RTS or decreasing RTS, but the magnitude of RTS cannot be determined.

(2) RTS measurement using BCC-DEA models. Banker et al. (1984) proposed the method to examine the intercept of the supporting hyperplane on the production possibility set (PPS) under a variable RTS assumption. This intercept corresponds to a dual variable regarding the convex constraint in BCC-DEA. The type of RTS can be determined by the sign of this intercept (positive, negative or zero). Note that the intercept interval need to be considered when the intercept is not unique. In addition, we can analyse the properties of DMU within a small neighbourhood to determine the RTS on this point. This type of research efforts includes Banker and Thrall (1992), Tone (1996), Golany and Yu (1997), Sueyoshi (1999), Cooper et al. (2000), Tone and Sahoo (2003). The main contribution of this type of research lies in providing a theoretical basis for not only the type of RTS but also the qualitative measurement of RTS.

(3) RTS measurement using FGL-DEA model and quantitative measurement of scale elasticity (SE). This type of research can be traced to the efforts of Färe and Grosskopf (1985) and Färe et al. (1983, 1985, 1994). They examined the scale efficiency to determine whether the DMU being evaluated achieves constant RTS.

Their approach identifies RTS through the ratios of a series of relative efficiencies obtained by different DEA models with radial measure, which has different constraints. There is no problem of multiple solutions for RTS, as occurs in the first two types of research. However, their approach is limited in determining the type of RTS. In this context, Førsund (1996) discussed the quantitative measurement of SE and RTS, which is extended further to mathematical characterisations of SE for both frontier and non-frontier units by Fukuyama (2000). Huang et al. (1997), Kerstens and Vanden Eeckaut (1998) and Read and Thanassoulis (2000) also investigated the quantitative measurement of SE in DEA models.

(4) RTS measurement using non-radial DEA models. There are a variety of DEA models. The most well-known DEA models are often referred to as "radial models", including the CCR-DEA model and BCC-DEA model with radial measure. The CCR-DEA model and BCC-DEA model make assumptions of constant RTS and variable RTS, respectively. The first three types of approaches are all based on radial models, which will miss slacks when evaluating DMUs. Therefore, scholars have proposed dozens of non-radial models (e.g., Zhu, 1996, 2001; Tone, 2001, 2002; Chen, 2003) to eliminate this problem, such as DEA models with Russell measures and additive models, among others. It is natural to explore the RTS measurement using non-radial models. For example, Banker et al. (2004) discussed the RTS measurement using an additive model and multiplicative model. Sueyoshi and Sekitani (2005) explored the RTS measurement using dynamic DEA whose production scheme includes a feedback process. Zarepisheh et al. (2010) discussed the RTS issue in multiplicative models, which is a single model in one stage and different with the two-stage method proposed by Banker et al. (2004). Lozano and Gutierrez (2011) analysed the RTS of 41 Spanish airports using DEA model with Russell measure. Khodabakhshi et al. (2010) discussed the RTS issue in vague DEA models. Sueyoshi and Sekitani (2007) theoretically explored the measurement of RTS using a non-radial model with a range-adjusted measure. A new linear programming RAM/RTS approach was proposed to address a simultaneous occurrence of multiple reference sets, multiple supporting hyperplanes and multiple projections. Soleimani-damaneh et al. (2006) explored the RTS measurement in FDH model. In fact, in the RTS measurement, the projections on the efficient frontier of DMUs within the production possibility set (PPS) are different between radial DEA models and non-radial DEA models, but this difference does not affect the RTS of DMUs on the efficient frontier.

Now we recall some basic facts on the classic RTS in DEA framework. Assume $DMU(X_0,Y_0) \in PPS$ $X_0 \in R_m^+, Y_0 \in R_s^+$, let $\beta(t) = \max\{\beta|((t+1)X_0,(\beta+1)Y_0) \in PPS\}$, where $t, \beta$ are input and output scaling factors, respectively. In the case that $\beta(t)$ is differentiable, the classical SE $e(X,Y)$ at any $(X,Y) = ((t+1)X_0,(\beta+1)Y_0) \in PPS$ is defined as the ratio of its marginal productivity (where it exists) to its average productivity, where marginal and average productivities are defined as $d\beta(t)/dt$ and $(\beta(t)+1)/(t+1)$ respectively (see, e.g., Podinovski et al., 2009). That is

$$e(X,Y) = \frac{d\beta(t)}{dt} \frac{t+1}{\beta(t)+1}$$

When $\beta(t) = t = 0$, we have the classical SE at $DMU(X_0,Y_0)$ as follows.

$$e(X_0,Y_0) = \frac{d\beta(t)}{dt}\bigg|_{t=0}$$

Banker (1984) firstly defined the classical RTS based on PPS in DEA framework. Podinovski et al. (2009), Podinovski and Førsund (2010) and Atici and Podinovski (2012) pointed out that the derivative in the above classical definition of SE (RTS) may not always exist, and thus they replaced the classical derivative by the directional derivatives, and defined left-hand and right-hand SE as follows.

**Definitions 5 and 6 (Left-hand and right-hand SE):** The left and right hand scale elasticity of $DMU(X_0, Y_0)$ are defined, respectively, as follows.

$$e^-(X_0,Y_0) = \frac{d\beta(t)}{dt}\bigg|_{t=0^-} \tag{5}$$

$$e^+(X_0,Y_0) = \frac{d\beta(t)}{dt}\bigg|_{t=0^+} \tag{6}$$

Then we can define:

(a) if $e^-(X_0,Y_0) > 1$ (or $e^+(X_0,Y_0) > 1$) holds, then increasing RTS prevails on the left-hand (or right-hand) side of this point;

(b) if $e^-(X_0,Y_0) = 1$ (or $e^+(X_0,Y_0) = 1$) holds, then constant RTS prevails on the left-hand (or right-hand) side of this point;

(c) if $e^-(X_0,Y_0) < 1$ (or $e^+(X_0,Y_0) < 1$) holds, then decreasing RTS prevails on the left-hand (or right-hand) side of this point.

**Remark 1:** *It is possible to combine left and right side scale elasticity to define an overall RTS as in Banker and Thrall (1992) and Podinovski and Førsund (2010). They defined that: (i) increasing RTS prevails at $DMU(X_0,Y_0)$ if and only if $e^-(X_0,Y_0) \geq e^+(X_0,Y_0) > 1$, (ii) constant RTS prevails at $DMU(X_0,Y_0)$ if and only if $e^-(X_0,Y_0) \geq 1 \geq e^+(X_0,Y_0)$, and (iii) decreasing RTS prevails at $DMU(X_0,Y_0)$ if and only if $1 > e^-(X_0,Y_0) \geq e^+(X_0,Y_0)$. In this paper we will use separately left and right directional RTS (see Section 3) to keep more information. If necessary, readers can combine them similarly as above.*

## 3 Directional SE and directional RTS

### 3.1 Definitions of directional SE and directional RTS for explicit production functions

Yang (2012) proposed the definitions of directional RTS and directional SE for explicit production functions. In this paper, we restate briefly and extend these definitions.

Let the input and output vectors be $X = (x_1, x_2, ..., x_m)^T$ and $Y = (y_1, y_2, ..., y_s)^T$, respectively. Assume we have a continuously differentiable mapping $F: R_+^{m+s} \rightarrow R^m$ given as follows:

$$F(X,Y) = 0, X \in R_+^m, Y \in R_+^s, \frac{\partial F(X,Y)}{\partial x_i} \leq 0, \frac{\partial F(X,Y)}{\partial y_j} \geq 0, r = 1,...,s, i = 1,...,m$$

where $F(X,Y)$ refers to the vector $(f_1(X,Y), f_2(X,Y), ..., f_s(X,Y))$, with $f_i$ being continuously differentiable.

**Remark 2:** *The above equation $F(X,Y) = 0$ will lead to an implicit mapping Y=g(X), which is the underlying production functions. However the above conditions are generally not enough to ensure the smoothness of the production functions. To this end, one needs to apply the implicit function theorems. For example, Krantz and Parks (2002) shows that if F(X,Y) satisfies the following conditions: i.e., F(X,Y) is a continuously differentiable function, and its Jacobian matrix is invertible, where the Jacobian matrix is defined as follows:*

$$[\bar{X}|\bar{Y}] = \begin{bmatrix} \frac{\partial f_1}{\partial x_1}(X,Y) & \cdots & \frac{\partial f_1}{\partial x_m}(X,Y) & \frac{\partial f_1}{\partial y_1}(X,Y) & \cdots & \frac{\partial f_1}{\partial y_s}(X,Y) \\ \vdots & \ddots & \vdots & \vdots & \ddots & \vdots \\ \frac{\partial f_s}{\partial x_1}(X,Y) & \cdots & \frac{\partial f_s}{\partial x_m}(X,Y) & \frac{\partial f_s}{\partial y_1}(X,Y) & \cdots & \frac{\partial f_s}{\partial y_s}(X,Y) \end{bmatrix}$$

where $\bar{X}$ is the matrix of partial derivatives in the variable $x_i(i=1,...,m)$, and $\bar{Y}$ is the matrix of partial derivatives in the variables $y_j(j=1,...,s)$, and then we can construct a mapping $g: R^m \to R^s$ whose graph $(X, g(X))$ is precisely the set of all $(X,Y)$ such that $F(X, g(X))=0$, and that $Y = g(X)$ is continuously differentiable.

In considering the idea of input and output change non-proportionally, we can express the input-output change as the following equation:

$$F(\beta_1 y_1,...,\beta_s y_s, t_1 x_1,...,t_m x_m) = 0 \tag{7}$$

where $(t_1,...,t_m)$ and $(\beta_1,...,\beta_s)$ represent the vectors of changes in input components of $X$, and the corresponding output components of $Y$, respectively.

Suppose we have

$$\begin{cases} \beta_r - 1 = \delta_r \beta + \varepsilon_r(t), r = 1,...,s \\ t_i - 1 = \omega_i t, i = 1,...,m \end{cases}$$

where $\beta$ and $t$ represent the amount of directional change of outputs and inputs, respectively. Parameters $\omega_i \geq 0, i=1,...,m$, $\sum_{i=1}^{m} \omega_i = m$ and $\delta_r \geq 0, r=1,...,s$, $\sum_{r=1}^{s} \delta_r = s$ are fixed numbers representing the directions of inputs and outputs (See **Remark 3**), respectively, and $\varepsilon_r(t), r=1,...,s$ are the higher order infinitesimals and satisfy $\lim_{t \to 0} \varepsilon_r(t) = 0$ and $\lim_{t \to 0} \varepsilon_r'(t) = 0$, respectively.

We define $\beta(t) = \max\{\beta : F(\beta_1 y_1,...,\beta_s y_s, t_1 x_1,...,t_m x_m) = 0\}$ and assume that the function $\beta(t)$ is defined in a very small neighbourhood of $t=0$. Firstly, let us examine the case where $\beta(t)$ is smooth. Similar to classical definition of SE (see, e.g., Podinovski et al., 2009), at any point let

$$(X,Y) = (diag\{1+\omega_1 t,...,1+\omega_m t\} X_0, \, diag\{\beta_1,...,\beta_s\} Y_0)$$

where $diag\{\cdot\}$ denotes the diagonal matrix.

We define its marginal productivity and average productivity as the outputs gained in direction $(\delta_1, \delta_2,...,\delta_s)^T$ by adding one unit of inputs in direction $(\omega_1, \omega_2,...,\omega_m)^T$ (denoted as $d\beta(t)/dt$) and the ratio of changes of outputs in direction $(\delta_1, \delta_2, ..., \delta_s)^T$ and inputs in direction $(\omega_1, \omega_2, ..., \omega_m)^T$ (denoted as $(1+t)/(1+\beta(t))$), respectively. Thus we can define its directional SE as the ratio of its marginal productivity to its average productivity. Thus we have

$$e(X,Y) = \frac{d\beta(t)}{dt} \frac{1+t}{1+\beta(t)}$$

In particular, at $(X_0, Y_0)$, where $t = \beta(t) = 0$, we have the following formula for directional SE:

$$e(X_0, Y_0) = \frac{d\beta(t)}{dt}\bigg|_{t=0} \quad (8)$$

The rationale behind formula (8) is as follows: If the quantity of the inputs is marginally increased by a factor $t > 0$ in direction $(\omega_1, \omega_2, ..., \omega_m)^T$, then the maximum quantity of the outputs possible in the technology will increase by a factor $t \times e(X_0, Y_0)$ in direction $(\delta_1, \delta_2, ..., \delta_s)^T$.

We can also differentiate (7) w.r.t. the input scaling factor $t$ and obtain:

$$\sum_{r=1}^{s} \frac{\partial F}{\partial y_r} y_r \frac{d\beta_r}{d\beta(t)} \frac{d\beta(t)}{dt} + \sum_{i=1}^{m} \frac{\partial F}{\partial x_i} x_i \frac{dt_i}{dt} = 0 \quad (9)$$

From Equation (9), we obtain:

$$e(X_0, Y_0) = \frac{d\beta(t)}{dt}\bigg|_{t=0} = -\sum_{i=1}^{m} \frac{\partial F}{\partial x_i} x_i \omega_i \bigg/ \sum_{r=1}^{s} \frac{\partial F}{\partial y_r} y_r \delta_r \bigg|_{(X_0,Y_0)} \quad (10)$$

Then, Equation (10) is the formula of directional SE at point $(X_0, Y_0)$ for the case where the production function is continuously differentiable in the directions of $(\omega_1, \omega_2, ..., \omega_m)^T$ and $(\delta_1, \delta_2, ..., \delta_s)^T$, where we assume $\omega_i \geq 0, i = 1, ..., m$, $\sum_{i=1}^{m} \omega_i = m$ and $\delta_r \geq 0, r = 1, ..., s$, $\sum_{r=1}^{s} \delta_r = s$.

Moreover, it is clear that in the smooth case for the diagonal direction (i.e., $\omega_i = 1, i = 1, ..., m; \delta_r = 1, r = 1, ..., s$), Equation (10) is as follows:

$$e(X_0, Y_0) = \frac{d\beta(t)}{dt}\bigg|_{t=0} = -\sum_{i=1}^{m} \frac{\partial F}{\partial x_i} x_i \bigg/ \sum_{r=1}^{s} \frac{\partial F}{\partial y_r} y_r \bigg|_{(X_0,Y_0)} \quad (11)$$

which is the same as the formula of traditional SE in economics (see, e.g., Førsund 1996).

However, since $\beta(t)$ is a maximum function, in general it may not be always differentiable even if *F(X,Y)* is smooth. In this case, applying the result in Borwein and Lewis (2006), there always exist the directional derivatives at *t=0*. Motivated by Equation (8) (as in Podinovski and Førsund (2010)), we can simply define the left and right-hand directional SE, respectively, as follows.

$$\text{Left: } e^-(X_0,Y_0) = \frac{d\beta(t)}{dt}\bigg|_{t=0^-} \quad (12)$$

$$\text{Right: } e^+(X_0,Y_0) = \frac{d\beta(t)}{dt}\bigg|_{t=0^+} \quad (13)$$

The traditional definitions of scale elasticity assume that the response function $\beta(t)$ is differentiable at $t=0$. However, Podinovski and Førsund (2010) and Atici and Podinovski (2012) pointed out that $\beta(t)$ is often not differentiable at $t=0$. Furthermore, for the case where the production is given by DEA, they demonstrated that its right-hand and left-hand derivatives always exist within the domain of $\beta(t)$. Readers are referred to these two papers for the details.

Therefore, we define:

(a) if $e^-(X_0,Y_0) > 1$ (or $e^+(X_0,Y_0) > 1$) holds, then increasing directional RTS prevails left-hand (or right-hand) of point $(X_0,Y_0)$ in the direction of $(\omega_1,\omega_2,...,\omega_m)^T$ and $(\delta_1,\delta_2,...,\delta_s)^T$;

(b) if $e^-(X_0,Y_0) = 1$ (or $e^+(X_0,Y_0) = 1$) holds, then constant directional RTS prevails left-hand (or right-hand) of point $(X_0,Y_0)$ in the direction of $(\omega_1,\omega_2,...,\omega_m)^T$ and $(\delta_1,\delta_2,...,\delta_s)^T$;

(c) if $e^-(X_0,Y_0) < 1$ (or $e^+(X_0,Y_0) < 1$) holds, then decreasing directional RTS prevails left-hand (or right-hand) of point $(X_0,Y_0)$ in the direction of $(\omega_1,\omega_2,...,\omega_m)^T$ and $(\delta_1,\delta_2,...,\delta_s)^T$.

***Remark 3:*** *Here we assume parameters $\sum_{i=1}^{m}\omega_i = m$ and $\sum_{r=1}^{s}\delta_r = s$. In fact we may assume $\sum_{i=1}^{m}\omega_i = A$ and $\sum_{r=1}^{s}\delta_r = B$, where A and B are arbitrary positive numbers respectively. In the given input and output direction, we can always have the directional SE, which depends on the vectors of input and output directions. Therefore, without loss of generality and for simplicity, we assume $\sum_{i=1}^{m}\omega_i = m$ and $\sum_{r=1}^{s}\delta_r = s$. In this case, for the diagonal direction, the directional SE is exactly the same as the formula of traditional SE in economics.*

## 3.2 Directional SE and RTS in DEA framework

In this subsection, we will introduce the definition of directional RTS into DEA framework, given the directional SE defined above for an explicit production function. Based on the definitions in Section 2 and Section 3.1, we propose the definitions of left-hand and right-hand directional SE on $DMU(X_0, Y_0)$ based on PPS as follows.

**Definitions 7 and 8 (Left-hand and right-hand directional SE in DEA):**

Assuming $DMU(X_0, Y_0) \in PPS$ and $X_0 \in R_m^+, Y_0 \in R_s^+$, we let

$$\beta(t) = \max\{\beta | (\Omega_t X_0, \Phi_\beta Y_0) \in PPS\}$$

where $\Omega_t = diag\{1+\omega_1 t,...,1+\omega_m t\}$ and $\Phi_\beta = diag\{1+\delta_1 \beta,...,1+\delta_s \beta\}$, $diag\{\cdot\}$ denotes the diagonal matrix. Vectors $(\omega_1,...,\omega_m)^T$ ($\omega_i \geq 0, i=1,...,m$) and $(\delta_1,...,\delta_s)^T$ ($\delta_r \geq 0, r=1,...,s$) represent inputs and outputs directions, respectively, and satisfy $\sum_{i=1}^{m}\omega_i = m; \sum_{r=1}^{s}\delta_r = s$ where $t, \beta$ are input and output scaling factors, respectively. The left-hand and right-hand directional SE on $DMU(X_0, Y_0)$ are as follows:

$$e^-(X_0, Y_0) = \frac{d\beta(t)}{dt}\bigg|_{t=0^-} \qquad (14)$$

$$e^+(X_0, Y_0) = \frac{d\beta(t)}{dt}\bigg|_{t=0^+} \qquad (15)$$

Then we have

(a) if $e^-(X_0, Y_0) > 1$ (or $e^+(X_0, Y_0) > 1$) holds, then increasing directional RTS prevails left-hand (or right-hand) of point $(X_0, Y_0)$ in the direction of $(\omega_1, \omega_2,...,\omega_m)$ and $(\delta_1, \delta_2,...,\delta_s)$;

(b) if $e^-(X_0, Y_0) = 1$ (or $e^+(X_0, Y_0) = 1$) holds, then constant directional RTS prevails left-hand (or right-hand) of point $(X_0, Y_0)$ in the direction of $(\omega_1, \omega_2,...,\omega_m)$ and $(\delta_1, \delta_2,...,\delta_s)$;

(c) if $e^-(X_0, Y_0) < 1$ (or $e^+(X_0, Y_0) < 1$) holds, then decreasing directional RTS prevails left-hand (or right-hand) of point $(X_0, Y_0)$ in the direction of $(\omega_1, \omega_2,...,\omega_m)$ and $(\delta_1, \delta_2,...,\delta_s)$.

It should be noted that there may exist some strongly efficient $(X, Y) \in PPS$ whose inputs cannot be further reduced in a direction of $(\omega_1, \omega_2,...,\omega_m)^T$ regardless of output direction $(\delta_1, \delta_2, ..., \delta_s)^T$. These strongly efficient DMUs are defined as the directional smallest scale size (DSSS, See **Definition 9**). On the contrary, because of the assumption of **free disposal**, the inputs of any strongly efficient $(X,Y) \in PPS$ could always be further expanded in a direction of $(\omega_1, \omega_2, ..., \omega_m)^T$. Therefore, we should address these two cases differently. Thus, we obtain the following **Definition 9** on the directional smallest scale size in the direction of $(\omega_1, \omega_2, ..., \omega_m)^T$ and $(\delta_1, \delta_2, ..., \delta_s)^T$.

**Definition 9:** The strongly efficient $(X_0, Y_0)$ is of the directional smallest scale size if and only if $(\Omega_t X_0, \Phi_\beta Y_0) \notin PPS$ for any $\beta$ and $t < 0$.

Banker and Thrall (1992) defined the extreme scale size (either the smallest scale size or the largest scale size) for weakly efficient units but without further elaboration on the largest scale size units. In our **Definition 9**, the directional smallest scale size is defined on the strongly efficient DMUs. Please note that the definition in our paper is compatible with the definition of extreme scale size in Banker and Thrall (1992) in the sense that their definition refers to the radial direction, the directional smallest scale size under **Definition 9** is also the smallest scale size defined in Banker and Thrall (1992).

## 4 Measurement of directional RTS

The DMUs are often divided into two categories in the measurement of RTS using DEA models. The two categories are referred to as the strongly efficient DMUs[1] on the efficient frontier and weakly efficient or inefficient DMUs. The RTS of weakly efficient or inefficient DMUs can be measured through their projections onto the strongly efficient frontier. This paper follows the above two categories and conducts RTS measurement based on the PPS produced by an input-based BCC-DEA model under the assumption of variable RTS and focuses on the RTS measurement of strongly efficient DMUs. The following Model (16) is input-based BCC-DEA with radial measure.

---

[1] Unless it is expressly stated, the efficient DMU refers to the strongly efficient DMU in this paper.

$$\min_{\theta_0, \lambda_j, s_i^-, s_r^+} \theta = \theta_0 - \varepsilon \left( \sum_{i=1}^{m} s_i^- + \sum_{r=1}^{s} s_r^+ \right)$$

$$s.t. \begin{cases} \sum_j \lambda_j x_{ij} + s_i^- = \theta_0 x_{i0}, i = 1,...,m \\ \sum_j \lambda_j y_{rj} - s_r^+ = y_{r0}, r = 1,...,s \\ \sum_j \lambda_j = 1, \lambda_j, s_i^-, s_r^+ \geq 0, r = 1,...,s; i = 1,...,m; j = 1,...,n \end{cases} \quad (16)$$

The dual form of Model (16) reads

$$\max_{u_r, v_i, \mu_0} \sum_{r=1}^{s} u_r y_{r0} + \mu_0$$

$$s.t. \begin{cases} \sum_{r=1}^{s} u_r y_{rj} - \sum_{i=1}^{m} v_i x_{ij} + \mu_0 \leq 0, j = 1,...,n \\ \sum_{i=1}^{m} v_i x_{i0} = 1 \\ u_r \geq 0, v_i \geq 0, r = 1,...,s, i = 1,...,m, \mu_0 \text{ free} \end{cases} \quad (17)$$

### 4.1 Finite difference method (FDM)

Golany and Yu (1997) used FDM to estimate RTS for each DMU by testing the existence of solutions in four regions defined in the neighbourhood of the analysed unit. They provided a procedure to determine the RTS to the "right" and "left" of the DMU being evaluated. Rosen et al. (1998) estimated the directional derivative of DMUs on strongly efficient frontiers using FDM. The basic idea of FDM is to examine the ratio of the amount of change of outputs $\beta$ on the efficient frontier in the specified direction caused by the increase (or decrease) in a small enough amount of inputs $t$ in the specified direction because RTS is a local property of DMUs.

### 4.1.1 Directional RTS measurement of strongly efficient DMUs

It is well known that the weakly or strongly efficient frontier of BCC-DEA is piecewise linear. Thus, we determine the directional RTS to the "right" and "left" of DMU being evaluated. Figure 3 shows the directional RTS to the "right" and "left" of the point E, which is on the strongly efficient frontier.

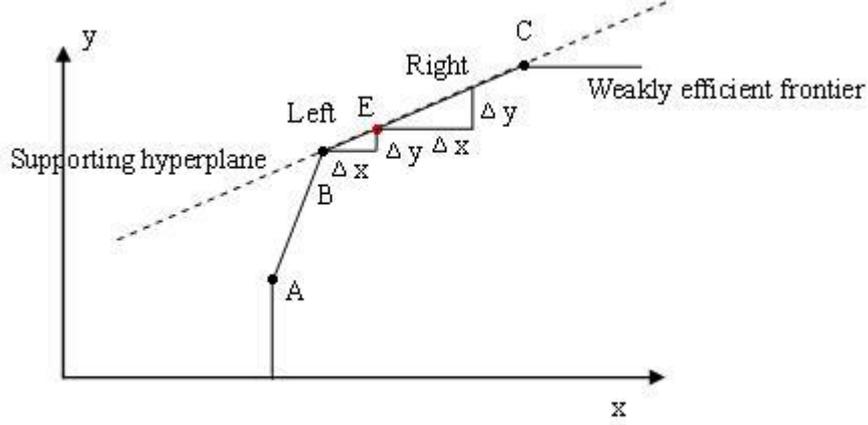

Figure 3: Directional RTS to the "right" and "left" of the point E

## 4.1.1.1 Directional RTS measurement to the "right" of strongly efficient DMUs

Based on the **Definition 8** and the FDM proposed by (Rosen et al., 1998; Golany and Yu, 1997), let $t_{right} > 0$. We have the following Model (18) to determine the right-hand directional RTS:

$$\max_{\beta, \lambda_j} \xi = \beta / t_{right}$$

$$s.t. \begin{cases} \sum_{j=1}^{n} \lambda_j x_{ij} \leq (1 + \omega_i t_{right}) x_{i0}, i = 1,...,m \\ \sum_{j=1}^{n} \lambda_j y_{rj} \geq (1 + \delta_r \beta) y_{r0}, r = 1,...,s \\ \sum_{j=1}^{n} \lambda_j = 1, \lambda_j \geq 0, j = 1,...,n \end{cases} \quad (18)$$

where $\delta_r \geq 0, r = 1,...,s$ and $\omega_i \geq 0, i = 1,...,m$ represent the direction factors of inputs and outputs, respectively, and satisfy $\sum_{r=1}^{s} \delta_r = s; \sum_{i=1}^{m} \omega_i = m$. Model (18) appears to be a nonlinear programming. However, as to be seen below, its objective is independent of $t_{right} > 0$ after $t_{right}$ is small enough. Thus, we will understand that $t_{right}$ is a small positive quantity, which represents the amount of directional change of inputs. Variable $\beta$ represents the amount of directional change of outputs. Then, it actually becomes a linear programming.

Let $\Omega_t^+ = diag\{1 + \omega_1 t_{right},...,1 + \omega_m t_{right}\}$, and $\Phi_{\beta^*}^+ = diag\{1 + \delta_1 \beta^*,...,1 + \delta_s \beta^*\}$, where $\beta^*$ is the optimal solution of Model (18). We have the following **Theorem 1**:

**Theorem 1.** There exists $t_0^+ > 0$ satisfying (1) when $t_{right} \in (0, t_0^+]$ and $(\Omega_t^+ X_0, \Phi_{\beta^*}^+ Y_0) \in PPS$, $(\Omega_t^+ X_0, \Phi_{\beta^*}^+ Y_0)$ is located on the weakly efficient frontier $EF_{weak}$ and (2) when $t_{right} \in (0, t_0^+]$, $(X_0, Y_0)$ and $(\Omega_t^+ X_0, \Phi_{\beta^*}^+ Y_0)$ have the same supporting hyperplane, which may be different for different $t_{right}$.

**Proof.** Please see Appendix A.

**Definition 10:** Let $\xi^*(X_0,Y_0)$ be the optimal objective of Model (18). Thus, we can determine the directional RTS to the "right" of $\text{DMU}(X_0,Y_0)$ as follows.

(a) if $\xi^*(X_0,Y_0)>1$ holds, then increasing directional RTS prevails in the direction of $(\omega_1,\omega_2,...,\omega_m)^T$ and $(\delta_1,\delta_2,...,\delta_s)^T$;

(b) if $\xi^*(X_0,Y_0)=1$ holds, then constant directional RTS prevails in the direction of $(\omega_1,\omega_2,...,\omega_m)^T$ and $(\delta_1,\delta_2,...,\delta_s)^T$;

(c) if $\xi^*(X_0,Y_0)<1$ holds, then decreasing directional RTS prevails in the direction of $(\omega_1,\omega_2,...,\omega_m)^T$ and $(\delta_1,\delta_2,...,\delta_s)^T$.

Here, we show that for very small positive $t_{right}$, the objective value in Model (18) is a constant for any given input and output direction.

**Theorem 2.** There exists a small enough quantity $t_0^+>0$ such that the optimal value of Model (18) is constant for all $t_{right}\in\left(0,t_0^+\right]$.

**Proof.** Please see Appendix A.

Now, we discuss how to select $t_{right}$ in practice. Consider the following Model (19):

$$\max_{\beta,\lambda_j} \beta$$
$$s.t.\begin{cases} \sum_{j=1}^n \lambda_j x_{ij} \leq (1+\omega_i t_{right})x_{i0}, i=1,...,m \\ \sum_{j=1}^n \lambda_j y_{rj} \geq (1+\delta_r \beta)y_{r0}, r=1,...,s \\ \sum_{j=1}^n \lambda_j = 1, \lambda_j \geq 0, j=1,...,n \end{cases} \quad (19)$$

Let $(\lambda_j^*,\beta^*)$ be optimal solutions of Model (19). Consider the following Model (20):

$$\max_{U,V,\mu_0} \varphi = U^T Y_0 + \mu_0$$
$$s.t.\begin{cases} U^T Y_j - V^T X_j + \mu_0 \leq 0, j=1,...,n \\ V^T X_0 = 1 \\ U^T \Phi_{\beta^*}^+ Y_0 - V^T \Omega_t^+ X_0 + \mu_0 = 0 \\ U \geq \mathbf{0}, V \geq \mathbf{0}, \mu_0 \text{ free} \end{cases} \quad (20)$$

From the proof of **Theorem 2**, if the optimal objective value of Model (20) satisfies $\varphi^*=1$, positive constant $t_{right}$ is small enough to guarantee that both $(X_0, Y_0)$ and $\left(\Omega_t^+ X_0, \Phi_{\beta^*}^+ Y_0\right)$ are located on the weakly efficient frontier $EF_{weak}$, and they have the same supporting hyperplane. Thus, in practice, we first select a small number $t_{right}$ in (19) and solve (20) to see if the optimal is the unit. If not, we will attempt smaller numbers. From the continuity, it will be the unit when $t_{right}$ is small enough.

**4.1.1.2 Directional RTS measurement to the "left" of strongly efficient DMUs**

First, we need to determine whether the strongly efficient $(X_0, Y_0)$ is of the directional smallest scale size. According to **Definition 9**, we consider the following Model (21):

$$\max_{\beta, \lambda_j, \eta} \eta$$
$$s.t. \begin{cases} \sum_{j=1}^n \lambda_j x_{ij} \leq (1-\omega_i \eta) x_{i0}, i=1,...,m \\ \sum_{j=1}^n \lambda_j y_{rj} \geq (1-\delta_r \beta) y_{r0}, r=1,...,s \\ \sum_{j=1}^n \lambda_j = 1, \lambda_j \geq 0, j=1,...,n \\ \eta \geq 0; \beta \text{ free} \end{cases} \quad (21)$$

**Theorem 3.** The optimal objective value $\eta^*$ of Model (21) is zero if and only if strongly efficient $(X_0, Y_0)$ is of the directional smallest scale size.

**Proof.** Please see Appendix A.

We first discuss the case in which strongly efficient $(X_0, Y_0)$ is not of the directional smallest scale size in the direction of $(\omega_1, \omega_2, ..., \omega_m)^T$ and $(\delta_1, \delta_2, ..., \delta_s)^T$. Based on the **Definition 7** and the FDM proposed by (Rosen et al., 1998; Golany and Yu, 1997), we let $t_{left}$ be a small positive constant and have the following Model (22) to determine the left-hand directional RTS:

$$\min_{\beta, \lambda_j} \psi = \beta/t_{left}$$
$$s.t. \begin{cases} \sum_{j=1}^n \lambda_j x_{ij} \leq (1-\omega_i t_{left}) x_{i0}, i=1,...,m \\ \sum_{j=1}^n \lambda_j y_{rj} \geq (1-\delta_r \beta) y_{r0}, r=1,...,s \\ \sum_{j=1}^n \lambda_j = 1, \lambda_j \geq 0, j=1,...,n \end{cases} \quad (22)$$

where $\delta_r \geq 0, r=1,...,s$ and $\omega_i \geq 0, i=1,...,m$ represent the direction factors of inputs and outputs, respectively, and satisfy $\sum_{r=1}^s \delta_r = s; \sum_{i=1}^m \omega_i = m$. Constant $t_{left}$ is a small positive quantity, which represents the amount of directional change of inputs. Variable $\beta$ represents the amount of directional change of outputs.

We let $\Omega_t^- = diag\{1-\omega_1 t_{left},...,1-\omega_m t_{left}\}$, $\Phi_{\beta^*}^- = diag\{1-\delta_1\beta^*,...,1-\delta_s\beta^*\}$, and $\beta^*$ is the optimal solution of Model (22). Thus, we have the following **Theorem 4**:

**Theorem 4.** There exists $t_0^- > 0$ satisfying (1) when $t_{left} \in (0, t_0^-]$ and $(\Omega_t^- X_0, \Phi_{\beta^*}^- Y_0) \in PPS$, $(\Omega_t^- X_0, \Phi_{\beta^*}^- Y_0)$ is located on the weakly efficient frontier $EF_{weak}$ and (2) when $t_{left} \in (0, t_0^-]$, $(X_0, Y_0)$ and $(\Omega_t^- X_0, \Phi_{\beta^*}^- Y_0)$ have the same supporting hyperplane, which may be different for different $t_{left}$.

**Proof.** Please see Appendix A.

**Definition 11:** We let $\psi^*(X_0, Y_0)$ be the optimal objective of Model (22). Accordingly, we can determine the directional RTS to the "left" of $DMU(X_0, Y_0)$ as follows:

(a) if $\psi^*(X_0, Y_0) > 1$ holds, then increasing directional RTS prevails in the direction of $(\omega_1, \omega_2, ..., \omega_m)^T$ and $(\delta_1, \delta_2, ..., \delta_s)^T$;

(b) if $\psi^*(X_0, Y_0) = 1$ holds, then constant directional RTS prevails in the direction of $(\omega_1, \omega_2, ..., \omega_m)^T$ and $(\delta_1, \delta_2, ..., \delta_s)^T$;

(c) if $\psi^*(X_0, Y_0) < 1$ holds, then decreasing directional RTS prevails in the direction of $(\omega_1, \omega_2, ..., \omega_m)^T$ and $(\delta_1, \delta_2, ..., \delta_s)^T$.

Next, we discuss how to choose $t_{left}$. Again, we first have

**Theorem 5.** There exists a small enough quantity $t_0^- > 0$ such that the optimal value of Model (22) is constant for all $t_{left} \in (0, t_0^-]$.

**Proof.** Please see Appendix A.

Now, again, consider the following Model (23):

$$\max_{\beta, \lambda_j} \beta$$

$$s.t. \begin{cases} \sum_{j=1}^{n} \lambda_j x_{ij} \leq (1-\omega_i t_{left}) x_{i0}, i=1,...,m \\ \sum_{j=1}^{n} \lambda_j y_{rj} \geq (1-\delta_r \beta) y_{r0}, r=1,...,s \\ \sum_{j=1}^{n} \lambda_j = 1, \lambda_j \geq 0, j=1,...,n \end{cases} \quad (23)$$

Let $(\lambda_j^*, \beta^*)$ be the optimal solutions of Model (23) and $\Omega_t^- = diag\{1-\omega_1 t_{left},...,1-\omega_m t_{left}\}$ and $\Phi_{\beta^*}^- = diag\{1-\delta_1\beta^*,...,1-\delta_s\beta^*\}$. Consider the

following Model (24):

$$\max_{U,V,\mu_0} \varphi = U^T Y_0 + \mu_0$$

$$s.t. \begin{cases} U^T Y_j - V^T X_j + \mu_0 \leq 0, j=1,...,n \\ V^T X_0 = 1 \\ U^T \Phi_t^- Y_0 - V^T \Omega_{\beta^*}^- X_0 + \mu_0 = 0 \\ U \geq \mathbf{0}, V \geq \mathbf{0}, \mu_0 \text{ free} \end{cases} \qquad (24)$$

Again, if the optimal objective value of Model (24) satisfies $\varphi^* = 1$, $t_{left} > 0$ is a small enough constant to guarantee that both $(X_0, Y_0)$ and $(\Omega_t^- X_0, \Phi_{\beta^*}^- Y_0)$ are located on the weakly efficient frontier $EF_{weak}$, and they have the same supporting hyperplane. Thus, we will select $t_{left} > 0$ similarly.

Now, we turn to the case that the strongly efficient $(X_0, Y_0)$ is of the directional smallest scale size in the direction of $(\omega_1, \omega_2, ..., \omega_m)^T$ and $(\delta_1, \delta_2, ..., \delta_s)^T$. In this case, we cannot find a feasible solution in Model (22) when $t_{left} > 0$ is a small positive constant. Thus, we provide the following **Definition 12** to address the left-hand directional RTS of $(X_0, Y_0)$:

**Definition 12:** If strongly efficient $(X_0, Y_0)$ is of the directional smallest scale size in the direction of $(\omega_1, \omega_2, ..., \omega_m)^T$ and $(\delta_1, \delta_2, ..., \delta_s)^T$, then increasing directional RTS prevails at the left-hand of $(X_0, Y_0)$.

Similar to Equation (A-5), we have

$$\psi^* = \frac{\beta^*}{t_{left}} = \frac{\sum_{i=1}^{m} v_i^* \omega_i x_{i0}}{\sum_{r=1}^{s} u_r^* \delta_r y_{r0}} \qquad (25)$$

where $U^* = (u_1^*, u_2^*, ..., u_s^*)^T$ and $V^* = (v_1^*, v_2^*, ..., v_m^*)^T$ are the optimal solutions of Model (24), and $(U^*, -V^*)$ is the normal vector of the supporting hyperplane on the $\text{DMU}(X_0, Y_0)$ and $\text{DMU}(\Omega_t^- X_0, \Phi_{\beta^*}^- Y_0)$.

### 4.1.1.3 A procedure for estimating directional RTS of strongly efficient DMUs

Based on the above analysis, we now propose a procedure for estimating directional RTS to the "right" and "left" of strongly efficient $\text{DMU}(X_0, Y_0)$ on the strongly efficient frontier $EF_{strong}$ as follows.

**Procedure 1.**

**Step 1: Determine the directional RTS to the "right" of $\mathrm{DMU}(X_0,Y_0)$**

**Step 1-1:** Choose a small enough quantity $t_{right}>0$, based on Model (19) and Model (20), to guarantee that both $(X_0,Y_0)$ and $(\Omega_t^+ X_0, \Phi_\beta^+ Y_0)$ are located on the weakly efficient frontier $EF_{weak}$, and they have the same supporting hyperplane.

**Step 1-2:** Solve Model (18) to determine the directional RTS to the "right" of $\mathrm{DMU}(X_0,Y_0)$:

(a) if $\xi^*(X_0,Y_0)>1$ holds, then increasing directional RTS prevails in the direction of $(\omega_1,\omega_2,...,\omega_m)^T$ and $(\delta_1,\delta_2,...,\delta_s)^T$;

(b) if $\xi^*(X_0,Y_0)=1$ holds, then constant directional RTS prevails in the direction of $(\omega_1,\omega_2,...,\omega_m)^T$ and $(\delta_1,\delta_2,...,\delta_s)^T$;

(c) if $\xi^*(X_0,Y_0)<1$ holds, then decreasing directional RTS prevails in the direction of $(\omega_1,\omega_2,...,\omega_m)^T$ and $(\delta_1,\delta_2,...,\delta_s)^T$.

**Step 2: Determine the directional RTS to the "left" of $\mathrm{DMU}(X_0,Y_0)$**

**Step 2-0:** Solve Model (21) to see whether its optimal objective value is zero. If so, $(X_0,Y_0)$ is of the directional smallest scale size. Otherwise, we have the following two steps to determine the left-hand directional RTS.

**Step 2-1:** Choose a small enough quantity $t_{left}$, based on Model (23) and Model (24), to guarantee that both $(X_0,Y_0)$ and $(\Omega_t^- X_0, \Phi_\beta^- Y_0)$ are located on the weakly efficient frontier $EF_{weak}$, and they have the same supporting hyperplane.

**Step 2-2:** Solve Model (22) to determine the directional RTS to the "left" of $\mathrm{DMU}(X_0,Y_0)$:

(a) if $\psi^*(X_0,Y_0)>1$ holds, then increasing directional RTS prevails in the direction of $(\omega_1,\omega_2,...,\omega_m)^T$ and $(\delta_1,\delta_2,...,\delta_s)^T$;

(b) if $\psi^*(X_0,Y_0)=1$ holds, then constant directional RTS prevails in the direction of $(\omega_1,\omega_2,...,\omega_m)^T$ and $(\delta_1,\delta_2,...,\delta_s)^T$;

(c) if $\psi^*(X_0,Y_0)<1$ holds, then decreasing directional RTS prevails in the direction of $(\omega_1,\omega_2,...,\omega_m)^T$ and $(\delta_1,\delta_2,...,\delta_s)^T$.

### 4.1.2 Directional RTS measurement of inefficient or weakly efficient DMUs

For estimating directional RTS to the "right" and "left" for inefficient or weakly efficient DMUs, we can perform the following two steps:

**Step 1:** First, we project the inefficient or weakly efficient DMUs onto the strongly efficient frontier $EF_{strong}$ using BCC-DEA with radial measure (See Model (16) for details), and the formula for projection is the following Equation (26):

$$\begin{cases} \tilde{x}_{i0} \leftarrow \theta_0^* x_{i0} - s_i^{-*}, i = 1,...,m \\ \tilde{y}_{r0} \leftarrow y_{r0} + s_r^{+*}, r = 1,...,s \end{cases} \quad (26)$$

**Step 2:** When we determine the projected points on the strongly efficient frontier $EF_{strong}$, we can estimate the directional RTS to the "right" and "left" for the inefficient or weakly efficient DMUs using **Procedure 1** in Section 4.1.1.3.

**Remark 4:** *Here we only use an input-based radial projection. However, different projections may generate different RTS (See, e.g., Sueyoshi and Sekitani, 2007). In this paper we concentrate on discussion of the directional RTS of strongly efficient DMUs instead of inefficient or weakly efficient ones so we only show one possible way to estimate their directional RTS by using the input-based radial projection.*

### 4.2 Upper and lower bounds method (ULBM)

In this subsection, we discuss how to estimate the directional RTS of the strongly efficient DMUs under the assumption of variable RTS from another viewpoint. According to Equations (A-5) and (25), we know that the following formulas hold when $t_{right} > 0$ and $t_{left} > 0$ are small enough positive constants:

$$\frac{\beta^*}{t_{right}} = \frac{\sum_{i=1}^m v_i^* \omega_i x_{i0}}{\sum_{r=1}^s u_r^* \delta_r y_{r0}} \quad \text{and} \quad \frac{\beta^*}{t_{left}} = \frac{\sum_{i=1}^m v_i^* \omega_i x_{i0}}{\sum_{r=1}^s u_r^* \delta_r y_{r0}} \quad (27)$$

Thus, we can use the following Model (28) to calculate the upper and lower bounds of the directional SE and then determine the type of directional RTS of $DMU(X_0, Y_0)$.

$$\bar{\rho}(\underline{\rho}) = \max_{u_r, v_i, \mu_0} \left( \min_{u_r, v_i, \mu_0} \right) \frac{\sum_{i=1}^m v_i \omega_i x_{i0}}{\sum_{r=1}^s u_r \delta_r y_{r0}}$$

$$s.t. \begin{cases} \sum_{r=1}^s u_r y_{rj} - \sum_{i=1}^m v_i x_{ij} + \mu_0 \leq 0, j = 1,...,n \\ \sum_{r=1}^s u_r y_{r0} - \sum_{i=1}^m v_i x_{i0} + \mu_0 = 0 \\ \sum_{i=1}^m v_i x_{i0} = 1 \\ u_r \geq 0, v_i \geq 0, r = 1,...,s, i = 1,...,m, \mu_0 \text{ free} \end{cases} \quad (28)$$

**Theorem 6.** Suppose that $(X_0, Y_0)$ is not of the directional smallest scale size. The upper and lower bounds ($\bar{\rho}(X_0, Y_0)$ and $\underline{\rho}(X_0, Y_0)$) in Model (28) are equal to the optimal objective value $\psi^*(X_0, Y_0)$ in Model (22) and the optimal objective value $\xi^*(X_0, Y_0)$ in Model (18).

**Proof.** Please see Appendix A.

**Theorem 7.** If the maximal optimal objective function $\bar{\rho}(X_0, Y_0)$ of Model (28) is unbounded ($+\infty$), the strongly efficient $(X_0, Y_0)$ is of the directional smallest scale size.

**Proof.** Please see Appendix A.

Therefore, we have the procedure similar to **Procedure 1** for estimating directional RTS.

Model (28) is a fractional programming that is difficult to solve, so we transform Model (28) into an equivalent linear programming through Charnes-Cooper transformation (Charnes et al., 1962). First, we rewrite Model (28) using the following Model (29):

$$\bar{\rho}(\underline{\rho}) = \max_{U,V,\mu_0} \left( \min_{U,V,\mu_0} \right) \frac{V^T W X_0}{U^T \Delta Y_0}$$

$$s.t. \begin{cases} U^T Y_j - V^T X_j + \mu_0 \leq 0, j=1,...,n \\ U^T Y_0 - V^T X_0 + \mu_0 = 0 \\ V^T X_0 = 1 \\ U \geq \mathbf{0}, V \geq \mathbf{0}, \mu_0 \text{ free} \end{cases} \quad (29)$$

where $U = (u_1, u_2, ..., u_s)^T$ and $V = (v_1, v_2, ..., v_m)^T$ are vectors of multipliers, and $\Delta = diag\{\delta_1, \delta_2, ..., \delta_s\}$ and $W = diag\{\omega_1, \omega_2, ..., \omega_m\}$ are matrixes of directions of inputs and outputs.

We let

$$\begin{cases} \tau = 1/U^T \Delta Y_0 \\ \tau V^T W = \Gamma^T \\ \tau U^T \Delta = \Lambda^T \end{cases} \quad (30)$$

We let $\omega_i \geq \varepsilon$, where $\varepsilon$ is a non-Archimedean construct to assure the inverse matrix of $W$ exists. In this case, we have the following Equation (31) from Equation (30):

$$\begin{cases} \tau = 1/U^T \Delta Y_0 \\ V^T = \dfrac{1}{\tau} \Gamma^T W^{-1} \\ U^T = \dfrac{1}{\tau} \Lambda^T \Delta^{-1} \end{cases} \tag{31}$$

Thus, Model (29) can be translated into the following Model (32):

$$\bar{\rho}(\underline{\rho}) = \max_{\Gamma,\Lambda,\tau,\mu_0} \left( \min_{\Gamma,\Lambda,\tau,\mu_0} \right) \Gamma^T X_0$$

$$s.t. \begin{cases} \Lambda^T \Delta^{-1} Y_j - \Gamma^T W^{-1} X_j + \tau \mu_0 \leq 0, j=1,\ldots,n \\ \Lambda^T \Delta^{-1} Y_0 - \Gamma^T W^{-1} X_0 + \tau \mu_0 = 0 \\ \Gamma^T W^{-1} X_0 = \tau \\ \Lambda^T Y_0 = 1 \\ \Gamma \geq \mathbf{0}, \Lambda \geq \mathbf{0}, \tau \geq 0, \mu_0 \text{ free} \end{cases} \tag{32}$$

We let $\tau \mu_0 = \mu_0'$, and Model (32) can be converted into the following linear programming:

$$\bar{\rho}(\underline{\rho}) = \max_{\Gamma,\Lambda,\tau,\mu_0'} \left( \min_{\Gamma,\Lambda,\tau,\mu_0'} \right) \Gamma^T X_0$$

$$s.t. \begin{cases} \Lambda^T \Delta^{-1} Y_j - \Gamma^T W^{-1} X_j + \mu_0' \leq 0, j=1,\ldots,n \\ \Lambda^T \Delta^{-1} Y_0 - \Gamma^T W^{-1} X_0 + \mu_0' = 0 \\ \Gamma^T W^{-1} X_0 = \tau \\ \Lambda^T Y_0 = 1 \\ \Gamma \geq \mathbf{0}, \Lambda \geq \mathbf{0}, \tau \geq 0, \mu_0' \text{ free} \end{cases} \tag{33}$$

Solving Model (33), we can obtain the optimal objective value of Model (28) or Model (29).

**5 A Case Study**

In this section, we conduct a case study to analyse the directional RTS of 16 basic research institutes in the Chinese Academy of Sciences (CAS) in 2010. Since the Pilot Project of Knowledge Innovation (KIPP) in 1998 at the CAS, institute evaluation has become increasingly important, and the requirements for the evaluation process have diversified. Since 2005, CAS headquarters has built up the Comprehensive Quality

Evaluation (CQE) system for institute evaluation in CAS. The results of evaluation are expressed as multi-dimensional feedback data and used as the tools to provide basis of comprehensive analysis and decision-making and to provide institutes with targeted evaluation information and diagnostic comments.

In the framework of CQE, multiple inputs and outputs of the basic research institutes of the CAS are monitored using several quantitative indicators. In this paper, we use the same index indicators as proposed in Liu et al. (2011) for 16 basic research institutes in CAS:

Inputs: (1) Staff denotes the number of full-time research staff, and (2) Res. Expen. denotes the total research expenditures;

Outputs: (1) SCI Pub. denotes the publications, including the international papers indexed by Science Citation Index; (2) High Pub. denotes high-quality papers published in top research journals; (3) Exter. Fund denotes the external research funding; (4) Grad. Enroll. denotes graduate students' enrolment.

Columns 2-7 of Table 2 show the detailed data of these indicators of 16 basic CAS research institutes in 2010.

Table 2: Detailed data of 16 basic CAS research institutes in 2010

| Institutes | Outputs | | | | Inputs | | Efficiency scores (Model 16) |
|---|---|---|---|---|---|---|---|
| | SCI Pub. (Number) | High Pub. (Number) | Grad. Enroll. (Number) | Exter. Fund. (RMB million) | Staff (FTE) | Res. Expen. (RMB million) | |
| $DMU_1$ | 436 | 133 | 184 | 31.5580 | 252 | 117.9450 | 1 |
| $DMU_2$ | 243 | 127 | 43 | 15.3041 | 37 | 29.4310 | 1 |
| $DMU_3$ | 164 | 70 | 89 | 33.8365 | 240 | 101.4250 | 0.6846 |
| $DMU_4$ | 810 | 276 | 247 | 183.8434 | 356 | 368.4830 | 1 |
| $DMU_5$ | 200 | 55 | 111 | 12.9342 | 310 | 195.8620 | 0.4538 |
| $DMU_6$ | 104 | 49 | 33 | 60.7366 | 201 | 188.8290 | 0.6179 |
| $DMU_7$ | 113 | 49 | 45 | 72.5368 | 157 | 131.3010 | 0.9597 |
| $DMU_8$ | 8 | 1 | 44 | 23.7015 | 236 | 77.4390 | 0.5300 |
| $DMU_9$ | 371 | 118 | 89 | 216.9885 | 805 | 396.9050 | 0.7777 |
| $DMU_{10}$ | 607 | 216 | 168 | 88.5561 | 886 | 411.5390 | 0.4541 |
| $DMU_{11}$ | 314 | 49 | 89 | 45.3597 | 623 | 221.4280 | 0.3575 |
| $DMU_{12}$ | 261 | 79 | 131 | 41.1156 | 560 | 264.3410 | 0.3702 |
| $DMU_{13}$ | 627 | 168 | 346 | 645.4150 | 1344 | 900.5090 | 1 |
| $DMU_{14}$ | 971 | 518 | 335 | 205.4528 | 508 | 344.3120 | 1 |
| $DMU_{15}$ | 395 | 180 | 117 | 90.0373 | 380 | 161.3310 | 0.8840 |
| DMU$_{16}$ | 229 | 138 | 62 | 32.6111 | 132 | 83.9720 | 0.6684 |

Data source: (1) Quantitative monitoring report of research institutes in CAS, 2011; (2) Statistical Yearbook of CAS, 2011.

Next, we analyse the directional RTS of 16 basic CAS research institutes in 2010.

**Step 1:** We determine the strongly efficient frontier $EF_{strong}$ and weakly efficient frontier $EF_{weak}$ using the input-based BCC-DEA model (16).

The last column of Table 2 shows the relative efficiencies of 16 basic CAS research institutes in 2010. From Table 2, we can see that DMU$_1$, DMU$_2$, DMU$_4$, DMU$_{13}$ and DMU$_{14}$ are efficient DMUs. Columns 2-6 of Table 3 show the projections of these DMUs on the strongly efficient frontier $EF_{strong}$ using the following formula:

$$\begin{cases} \tilde{x}_{i0} \leftarrow \theta_0^* x_{i0} - s_i^{-*}, i=1,...,m \\ \tilde{y}_{r0} \leftarrow y_{r0} + s_r^{+*}, r=1,...,s \end{cases} \quad (34)$$

**Step 2:** From Models (18), (21), (22) and (33), we can obtain the directional RTS to the "right" and "left" of each DMU. For comparison purposes, we set the direction of outputs as $\delta_1 = \delta_2 = \delta_3 = \delta_4 = 1$ and the direction of inputs as $\omega_1 = 0.5, \omega_2 = 1.5$, $\omega_1 = 1.5, \omega_2 = 0.5$ and $\omega_1 = 1, \omega_2 = 1$ (traditional RTS), respectively. We use Model (21) to determine whether DMUs are of the directional smallest scale size. We find that DMU$_2$ is of the directional smallest scale size in the above three input directions. Therefore, we let $t_{right} = t_{left} = 1E^{-6}$, which can pass the tests of Equation (19) - (20) and Equation (23) - (24) under the above three directions. Columns 8-10 of Table 3 show the directional RTS of each DMU in the case of $\omega_1 = 0.5, \omega_2 = 1.5$ (Case 1), $\omega_1 = 1.5, \omega_2 = 0.5$ (Case 2) and the traditional case $\omega_1 = 1, \omega_2 = 1$ (Case 3), respectively.

In Table 3, I, C and D denote increasing, constant and decreasing directional RTS, respectively. Based on the results of analysis, we find that the directional RTS of DMUs may change due to different directions of inputs, and directional RTS is different than traditional RTS. Taking DMU$_{16}$ as an example, we can see that (1) decreasing directional RTS prevails to both the "left" and the "right" of this DMU in the direction of $\omega_1 = 1, \omega_2 = 1$ and $\delta_1 = \delta_2 = \delta_3 = \delta_4 = 1$; (2) increasing directional RTS prevails to the "left" and decreasing directional RTS prevails to the "right" of this DMU in the direction of $\omega_1 = 0.5, \omega_2 = 1.5$ and $\delta_1 = \delta_2 = \delta_3 = \delta_4 = 1$; and (3) constant directional RTS prevails to the "left" and decreasing directional RTS prevails to the "right" of this DMU in the direction of $\omega_1 = 1.5, \omega_2 = 0.5$ and $\delta_1 = \delta_2 = \delta_3 = \delta_4 = 1$. See Table 3 for details.

Table 3: The projections of these DMUs on the strongly efficient frontier and directional RTS in three cases

| Institutes | Projections (Outputs) | | | | Projections (Inputs) | | Case 1 (Left /Right) | Case 2 (Left /Right) | Case 3 (Left /Right) |
|---|---|---|---|---|---|---|---|---|---|
| | SCI Pub. (Number) | High Pub. (Number) | Grad. Enroll. (Number) | Exter. Fund. (RMB million) | Staff (FTE) | Res. Expen. (RMB million) | | | |
| $DMU_1$ | 436 | 133 | 184 | 31.5580 | 252 | 117.9450 | I/D | I/D | I/D |
| $DMU_2$ | 243 | 127 | 43 | 15.3041 | 37 | 29.4310 | I(DSSS)/D | I(DSSS)/D | I(DSSS)/D |
| $DMU_3$ | 333.7205 | 160.9626 | 89 | 33.8365 | 109.3189 | 69.4315 | I/D | D/D | D/D |
| $DMU_4$ | 810 | 276 | 247 | 183.8434 | 356 | 368.4830 | I/D | I/D | I/D |
| $DMU_5$ | 336.0780 | 129.8936 | 111 | 23.1429 | 140.6880 | 72.1186 | I/D | I/D | I/D |
| $DMU_6$ | 377.6452 | 161.7549 | 93.1720 | 60.7366 | 124.1866 | 116.6665 | D/D | I/D | I/D |
| $DMU_7$ | 354.3239 | 153.4511 | 90.7659 | 72.5368 | 150.6740 | 126.0097 | D/D | I/D | I/D |
| $DMU_8$ | 248.1175 | 127.5464 | 47.0380 | 23.7015 | 54.4182 | 41.0397 | I/D | I/D | I/D |
| $DMU_9$ | 371 | 143.2717 | 141.6513 | 216.9885 | 455.9783 | 308.6761 | I/D | I/D | C/D |
| $DMU_{10}$ | 607 | 322.5000 | 189 | 110.3785 | 272.5000 | 186.8715 | I/D | I/D | D/D |
| $DMU_{11}$ | 357.5985 | 188.5111 | 89.0740 | 45.3597 | 111.2927 | 791.0020 | D/D | I/I | D/D |
| $DMU_{12}$ | 396.3237 | 168.6467 | 131 | 41.1156 | 173.7626 | 97.8468 | I/D | I/D | D/D |
| $DMU_{13}$ | 627 | 168 | 346 | 645.4150 | 1344 | 900.5090 | I/D | I/D | I/D |
| $DMU_{14}$ | 971 | 518 | 335 | 205.4528 | 508 | 344.3120 | I/D | I/D | C/D |
| $DMU_{15}$ | 404.7128 | 203.7191 | 117 | 90.0373 | 206.5457 | 142.6131 | I/D | I/D | D/D |
| $DMU_{16}$ | 286.1352 | 148.2834 | 62 | 32.6111 | 76.9688 | 56.1180 | I/D | C/D | D/D |

Note: DSSS denotes that $DMU_2$ is of the directional smallest scale size, and increasing directional RTS prevails to the "left" of $DMU_2$ according to **Definition 12**.

**Step 3:** We take strongly efficient $DMU_1$ and $DMU_5$ as examples to analyse their directional RTS in multiple input directions. As before, we set the directions of outputs as $\delta_1 = \delta_2 = \delta_3 = \delta_4 = 1$. We apply Model (21) to these two DMUs and find that the optimal objective values are all larger than zero in different directions of inputs. We let $t_{right} = t_{left} = 1E^{-6}$, which can pass the tests of Equation (19) - (20) and Equation (23) - (24) in different directions of inputs. Thus, we can obtain the directional SE and RTS of $DMU_1$ and $DMU_5$ in different directions of inputs using FDM method. Additionally, we use Model (33) to obtain the upper and lower bounds of objective function value. See Table 4 for details.

Table 4: The directional RTS of $DMU_1$ and $DMU_5$ in different input directions

| DMU | $\omega_1$ | $\omega_2$ | $\xi^*$ (Right) | $\psi^*$ (Left) | $\underline{\rho}$ (Lower | $\overline{\rho}$ (Upper | Directional RTS(Right) | Directional RTS(Left) |
|---|---|---|---|---|---|---|---|---|

|   |   |   |   |   | bound) |   | bound) |   |   |
|---|---|---|---|---|---|---|---|---|---|
| | 0.1 | 1.9 | 0.08 | 1.94 | 0.08 | 1.94 | Decreasing | Increasing |
| | 0.2 | 1.8 | 0.17 | 1.84 | 0.17 | 1.84 | Decreasing | Increasing |
| | 0.3 | 1.7 | 0.25 | 1.74 | 0.25 | 1.74 | Decreasing | Increasing |
| | 0.4 | 1.6 | 0.33 | 1.63 | 0.33 | 1.63 | Decreasing | Increasing |
| | 0.5 | 1.5 | 0.41 | 1.53 | 0.41 | 1.53 | Decreasing | Increasing |
| | 0.6 | 1.4 | 0.49 | 1.43 | 0.49 | 1.43 | Decreasing | Increasing |
| | 0.7 | 1.3 | 0.56 | 1.33 | 0.56 | 1.33 | Decreasing | Increasing |
| | 0.8 | 1.2 | 0.51 | 1.23 | 0.51 | 1.23 | Decreasing | Increasing |
| | 0.9 | 1.1 | 0.47 | 1.12 | 0.47 | 1.12 | Decreasing | Increasing |
| $DMU_1$ | 1 | 1 | 0.43 | 1.02 | 0.43 | 1.02 | Decreasing | Increasing |
| | 1.1 | 0.9 | 0.38 | 1 | 0.38 | 1 | Decreasing | Constant |
| | 1.2 | 0.8 | 0.34 | 1.08 | 0.34 | 1.08 | Decreasing | Increasing |
| | 1.3 | 0.7 | 0.30 | 1.17 | 0.30 | 1.17 | Decreasing | Increasing |
| | 1.4 | 0.6 | 0.26 | 1.26 | 0.26 | 1.26 | Decreasing | Increasing |
| | 1.5 | 0.5 | 0.21 | 1.35 | 0.21 | 1.35 | Decreasing | Increasing |
| | 1.6 | 0.4 | 0.17 | 1.44 | 0.17 | 1.44 | Decreasing | Increasing |
| | 1.7 | 0.3 | 0.13 | 1.53 | 0.13 | 1.53 | Decreasing | Increasing |
| | 1.8 | 0.2 | 0.09 | 1.62 | 0.09 | 1.62 | Decreasing | Increasing |
| | 1.9 | 0.1 | 0.04 | 1.71 | 0.04 | 1.71 | Decreasing | Increasing |
| | 0.1 | 1.9 | 0.08 | 1.97 | 0.08 | 1.97 | Decreasing | Increasing |
| | 0.2 | 1.8 | 0.16 | 1.86 | 0.16 | 1.86 | Decreasing | Increasing |
| | 0.3 | 1.7 | 0.24 | 1.76 | 0.24 | 1.76 | Decreasing | Increasing |
| | 0.4 | 1.6 | 0.32 | 1.66 | 0.32 | 1.66 | Decreasing | Increasing |
| | 0.5 | 1.5 | 0.40 | 1.55 | 0.40 | 1.55 | Decreasing | Increasing |
| | 0.6 | 1.4 | 0.48 | 1.45 | 0.48 | 1.45 | Decreasing | Increasing |
| | 0.7 | 1.3 | 0.55 | 1.35 | 0.55 | 1.35 | Decreasing | Increasing |
| | 0.8 | 1.2 | 0.60 | 1.24 | 0.60 | 1.24 | Decreasing | Increasing |
| | 0.9 | 1.1 | 0.55 | 1.14 | 0.55 | 1.14 | Decreasing | Increasing |
| $DMU_5$ | 1 | 1 | 0.50 | 1.03 | 0.50 | 1.03 | Decreasing | Increasing |
| | 1.1 | 0.9 | 0.45 | 0.93 | 0.45 | 0.93 | Decreasing | Decreasing |
| | 1.2 | 0.8 | 0.40 | 1 | 0.40 | 1 | Decreasing | Constant |
| | 1.3 | 0.7 | 0.35 | 1.08 | 0.35 | 1.08 | Decreasing | Increasing |
| | 1.4 | 0.6 | 0.30 | 1.16 | 0.30 | 1.16 | Decreasing | Increasing |
| | 1.5 | 0.5 | 0.25 | 1.25 | 0.25 | 1.25 | Decreasing | Increasing |
| | 1.6 | 0.4 | 0.20 | 1.33 | 0.20 | 1.33 | Decreasing | Increasing |
| | 1.7 | 0.3 | 0.15 | 1.41 | 0.15 | 1.41 | Decreasing | Increasing |
| | 1.8 | 0.2 | 0.10 | 1.50 | 0.10 | 1.50 | Decreasing | Increasing |
| | 1.9 | 0.1 | 0.50 | 1.58 | 0.50 | 1.58 | Decreasing | Increasing |

From Table 4, we can see $\underline{\rho} = \xi^*$ and $\bar{\rho} = \psi^*$. Based on the above analysis, we have the following summaries, see Figures 4 to 7 for details.

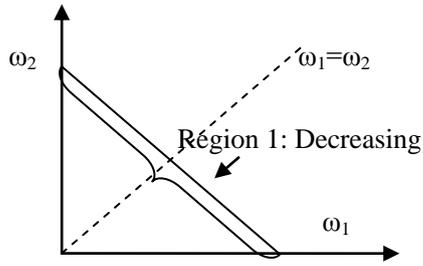
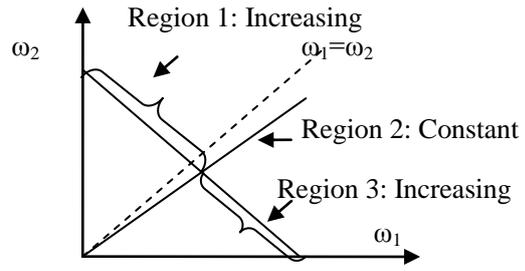

Figure 4: Directional RTS to the "right" of $DMU_1$    Figure 5: Directional RTS to the "right" of $DMU_5$

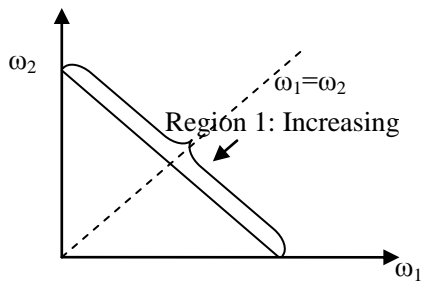
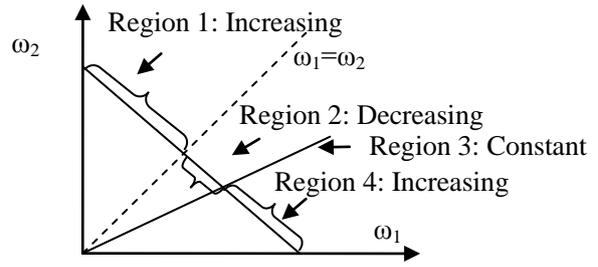

Figure 6: Directional RTS to the "left" of $DMU_1$    Figure 7: Directional RTS to the "left" of $DMU_5$

**(a) The directional RTS to the "right" of $DMU_1$ and $DMU_5$**

(a-1) For $DMU_1$, on the basis of existing inputs, if Staff and Res. Expen. increase in any proportion (under Pareto preference), decreasing directional RTS prevails on $DMU_1$, i.e., $DMU_1$ locates to the region with decreasing directional RTS in any direction of input increase. See Figure 4.

(a-2) For $DMU_5$, on the basis of existing inputs, if Staff and Res. Expen. increase in radial proportion, increasing directional RTS prevails on $DMU_5$. If the proportion of Staff and Res. Expen. increase locates in Region 1 and Region 3, increasing directional RTS prevails. Otherwise, if the proportion of inputs increase locates in Region 2 in Figure 5, constant directional RTS prevails. See Figure 5.

**(b) The directional RTS to the "left" of $DMU_1$ and $DMU_5$**

(b-1) For $DMU_1$, on the basis of existing inputs, if Staff and Res. Expen. decrease in any proportion (under Pareto preference), increasing directional RTS prevails on $DMU_1$. See Figure 6.

(b-2) For $DMU_5$, on the basis of existing inputs, if Staff and Res. Expen. decrease in radial proportion, increasing directional RTS prevails. If the proportion of Staff and Res. Expen. locates in Region 1 and Region 4 in Figure 7, increasing directional RTS prevails. If the proportion of inputs decrease locates in Region 2, decreasing directional RTS prevails. If the proportion of inputs decrease locates in Region 3, constant directional RTS prevails. See Figure 7.

**6 Conclusions**

In research institutions, given the complexity of research activities, changes of various types of inputs or outputs are often not proportional. Therefore, the existing definition of RTS in the framework of the DEA method may not meet the needs for estimation of the RTS of research institutions with multiple inputs and outputs using the DEA method. This work extends the definition of RTS in the DEA framework, proposes the definition of directional RTS in DEA and estimates the directional RTS of research institutions using DEA models. The traditional RTS is a special case of directional RTS when changes of input-output are in the diagonal direction. The directional RTS can be used for analysing those input change directions that are suitable for a particular DMU and thus useful for decision-makers (DMs) to decide rational combination of resources.

**Acknowledgements.** We would like to acknowledge the support of the National Natural Science Foundation of China (No. 71201158) and German Academic Exchange Service (DAAD, No. A1394033). We wish to express our sincere thanks to the referees and editors, whose suggestions greatly improved our paper's quality.

**Appendix A.**

   **Proof of Theorem 1.** (1) When $t_0^+ \geq t_{right} > 0$ and $\left(\Omega_t^+ X_0, \Phi_{\beta^*}^+ Y_0\right) \in PPS$, $\left(\Omega_t^+ X_0, \Phi_{\beta^*}^+ Y_0\right)$ is located on the weakly efficient frontier $EF_{weak}$. Assuming $\left(\Omega_t^+ X_0, \Phi_{\beta^*}^+ Y_0\right)$ is not located on the weakly efficient frontier $EF_{weak}$, we explore the following Model (A-1):

$$\max \theta$$
$$s.t. \begin{cases} \sum_{j=1}^n \lambda_j x_{ij} \leq \left(1 + \omega_i t_{right}\right) x_{i0}, i = 1, \ldots, m \\ \sum_{j=1}^n \lambda_j y_{rj} \geq \theta \left(1 + \delta_r \beta^*\right) y_{r0}, r = 1, \ldots, s \\ \sum_{j=1}^n \lambda_j = 1, \lambda_j \geq 0, j = 1, \ldots, n \end{cases} \quad (A-1)$$

Let $\theta^*$ be the optimal objective of the above model. Because $\left(\Omega_t^+ X_0, \Phi_{\beta^*}^+ Y_0\right)$ is not located on the weakly efficient frontier $EF_{weak}$, we have $\theta^* > 1$. Therefore, we have $\left(\Omega_t^+ X_0, \Theta^* \Phi_{\beta^*}^+ Y_0\right) \in PPS$, where

$$\Theta^* = diag\underbrace{\left\{\theta^*, \ldots, \theta^*\right\}}_{s}.$$

In addition, because $\theta^* > 1$, we have

$$\theta^* \left(1 + \delta_r \beta^*\right) y_{r0} \geq \left(1 + \delta_r \left(\theta^* \beta^*\right)\right) y_{r0}, r = 1, \ldots, s$$

Hence, $\theta^* \beta^*$ is a feasible solution for Model (18) and $\theta^* \beta^* > \beta^*$, which contradicts the fact that $\beta^*$ is the optimal solution of Model (18).

(2) When $t_{right} \to 0$, $\left(\Omega_t^+ X_0, \Phi_{\beta^*}^+ Y_0\right)$ converges to $(X_0, Y_0)$. Thus, there exists a small quantity $t_{right}$ satisfying that both $(X_0, Y_0)$ and $\left(\Omega_t^+ X_0, \Phi_{\beta^*}^+ Y_0\right)$ obviously having the same supporting hyperplane. This supporting hyperplane may be different for different $t_{right}$. **Q.E.D.**

   **Proof of Theorem 2.** We let $\left(\lambda_j^*, \beta^*\right)$ be the optimal solutions of Model (18). Because $t_{right} > 0$ is a small enough quantity, we know that both $(X_0, Y_0)$ and $\left(\Omega_t^+ X_0, \Phi_{\beta^*}^+ Y_0\right)$ are located on the weakly efficient frontier $EF_{weak}$, and they have the same supporting hyperplane, where $\Omega_t^+ = diag\{1 + \omega_1 t_{right}, \ldots, 1 + \omega_m t_{right}\}$ and $\Phi_{\beta^*}^+ = diag\{1 + \delta_1 \beta^*, \ldots, 1 + \delta_s \beta^*\}$. In this case, we know that the optimal objective value of Model (20) satisfies $\varphi^* = 1$. Thus, we have

$$\begin{cases} U^{T*} Y_0 - V^{T*} X_0 + \mu_0^* = 0 \\ U^{T*} \Phi_{\beta^*}^+ Y_0 - V^{T*} \Omega_t^+ X_0 + \mu_0^* = 0 \end{cases} \quad (A-2)$$

where $U^* = (u_1^*, u_2^*, ..., u_s^*)^T$ and $V^* = (v_1^*, v_2^*, ..., v_m^*)^T$ are the optimal solutions of Model (20).

According to Equation (A-2), we know

$$\sum_{r=1}^{s} u_r^* \delta_r \beta^* y_{r0} - \sum_{i=1}^{m} v_i^* \omega_i t_{right} x_{i0} = 0 \qquad (A-3)$$

From Equation (A-3), we have

$$\frac{\beta^*}{t_{right}} = \frac{\sum_{i=1}^{m} v_i^* \omega_i x_{i0}}{\sum_{r=1}^{s} u_r^* \delta_r y_{r0}} \qquad (A-4)$$

Therefore, when $t_{right} > 0$ is a small enough quantity, we can obtain the optimal objective value of Model (18) as follows:

$$\xi^* = \frac{\beta^*}{t_{right}} = \frac{\sum_{i=1}^{m} v_i^* \omega_i x_{i0}}{\sum_{r=1}^{s} u_r^* \delta_r y_{r0}} \qquad (A-5)$$

where $U^* = (u_1^*, u_2^*, ..., u_s^*)^T$ and $V^* = (v_1^*, v_2^*, ..., v_m^*)^T$ are the optimal solutions of Model (20), and $(U^*, -V^*)$ is the normal vector of the supporting hyperplane on the DMU$(X_0, Y_0)$ and DMU$(\Omega_t^+ X_0, \Phi_{\beta^*}^+ Y_0)$.

As the limit of Equation (15) always exists, we know when $0 < t_{right} \leq t_0^+$ and $t_{right} \to 0^+$, both $(X_0, Y_0)$ and $(\Omega_t^+ X_0, \Phi_{\beta^*}^+ Y_0)$ are located on the same "Face" of the weakly efficient frontier, or the value of the Equation (A-5) remains unchanged. Otherwise, the limit of Equation (15) does not exist. Thus, from Equation (A-5), we know that the optimal objective value of Model (18) is constant with respect to $t_{right}$ when $t_{right} > 0$ is a small enough quantity. **Q.E.D.**

**Proof of Theorem 3.** According to **Definition 9**, we can easily see that Theorem 3 holds. **Q.E.D.**

**Proof of Theorem 4.** The proof is similar to that of **Theorem 1** and omitted here.

**Proof of Theorem 5.** The proof is similar to that of **Theorem 2** and omitted here.

**Proof of Theorem 6.** We first discuss the equality between the optimal objective value $\xi^*(X_0, Y_0)$ and the lower bound $\underline{\rho}(X_0, Y_0)$.

From Equation (A-5), we know that the directional RTS to the "right" of DMU$(X_0, Y_0)$ reads

$$\xi^*(X_0, Y_0) = \frac{\beta^*}{t_{right}} = \frac{\sum_{i=1}^m v_i^* \omega_i x_{i0}}{\sum_{r=1}^s u_r^* \delta_r y_{r0}} \tag{A-6}$$

where $t_{right} > 0$ is a small enough quantity and $\{U^* = (u_1^*, u_2^*, ..., u_s^*)^T, V^* = (v_1^*, v_2^*, ..., v_m^*)^T, \mu_0^*\}$ is the optimal solution of Model (20).

(1) If $\xi^*(X_0, Y_0) < \underline{\rho}(X_0, Y_0)$, the optimal solution of Model (20) $\{U^* = (u_1^*, u_2^*, ..., u_s^*)^T, V^* = (v_1^*, v_2^*, ..., v_m^*)^T, \mu_0^*\}$ satisfies

$$\begin{cases} \sum_{r=1}^s u_r^* y_{rj} - \sum_{i=1}^m v_i^* x_{ij} + \mu_0^* \leq 0, j = 1,...,n \\ \sum_{r=1}^s u_r^* y_{r0} - \sum_{i=1}^m v_i^* x_{i0} + \mu_0^* = 0 \\ \sum_{i=1}^m v_i^* x_{i0} = 1 \end{cases} \tag{A-7}$$

which contradicts the fact that $\underline{\rho}(X_0, Y_0)$ is the lower bound in Model (28).

(2) If $\xi^*(X_0, Y_0) > \underline{\rho}(X_0, Y_0)$, we can deduce the following formula (A-8) from Equation (A-5) and Equation (27):

$$\xi^*(X_0, Y_0) = \frac{\beta^*}{t_{right}} > \underline{\rho}(X_0, Y_0) = \frac{\beta_{\underline{\rho}}^*}{t_{right}} \tag{A-8}$$

As $t_{right} > 0$, we have $\beta^* > \beta_{\underline{\rho}}^*$.

From Model (28), we know that $DMU(X_0, Y_0)$ satisfies

$$U^{*T} Y_0 - V^{*T} X_0 + \mu_0^* = 0 \tag{A-9}$$

$$U_{\underline{\rho}}^{*T} Y_0 - V_{\underline{\rho}}^{*T} X_0 + \mu_{\underline{\rho}0}^* = 0 \tag{A-10}$$

where $(U^*, -V^*)$ is the normal vector of a certain "Face" of weakly efficient frontier on the $DMU(X_0, Y_0)$, and $(U_{\underline{\rho}}^*, -V_{\underline{\rho}}^*)$ is the normal vector of a supporting hyperplane on $DMU(X_0, Y_0)$.

From (A-6) and (A-9), we have

$$U^{T*} \Phi_{\beta^*}^+ Y_0 - V^{T*} \Omega_t^+ X_0 + \mu_0 = 0 \tag{A-11}$$

where $\Omega_t^+ = diag\{1 + \omega_1 t_{right}, ..., 1 + \omega_m t_{right}\}$ and $\Phi_{\beta^*}^+ = diag\{1 + \delta_1 \beta^*, ..., 1 + \delta_s \beta^*\}$.

We pick up a point $\left(\Omega_t^+ X_0, \Phi_{\beta_\rho^*}^+ Y_0\right)$ on the supporting hyperplane with the normal vector $\left(U_\rho^*, -V_\rho^*\right)$ on the DMU $(X_0, Y_0)$. Thus, we obtain

$$U_\rho^{T*}\Phi_{\beta_\rho^*}^+ Y_0 - V_\rho^{T*}\Omega_t^+ X_0 + \mu_{\rho 0}^* = 0 \qquad (A\text{-}12)$$

where $\Phi_{\beta_\rho^*}^+ = diag\{1+\delta_1\beta_\rho^*, ..., 1+\delta_s\beta_\rho^*\}$. As DMU $(X_0, Y_0)$ is strongly efficient, there exists at least one set $(U_\rho^*, V_\rho^*) > \mathbf{0}$ that satisfies $U_\rho^{T*}Y_0 - V_\rho^{T*}X_0 + \mu_{\rho 0}^* = 0$. In this context, we can obtain the point $\left(\Omega_t^+ X_0, \Phi_{\beta_\rho^*}^+ Y_0\right)$.

As $\beta^* > \beta_\rho^*$, we can obtain the following formula from Equation (A-12):

$$U_\rho^{T*}\Phi_{\beta^*}^+ Y_0 - V_\rho^{T*}\Omega_t^+ X_0 + \mu_{\rho 0}^* > 0 \qquad (A\text{-}13)$$

We know that $\left(\Omega_t^+ X_0, \Phi_{\beta^*}^+ Y_0\right)$ is on the weakly efficient frontier $EF_{weak}$, and this fact contradicts the supporting hyperplane $U_\rho^{T*}Y - V_\rho^{T*}X + \mu_{\rho 0}^* = 0$ in Model (28).

Similarly, we can prove that the optimal objective value of Model (22) is equal to the upper bound $\bar{\rho}(X_0, Y_0)$ in Model (28). **Q.E.D.**

**Proof of Theorem 7.** We suppose that the maximal optimal objective function $\bar{\rho}(X_0, Y_0)$ of Model (28) is unbounded, but strongly efficient $(X_0, Y_0)$ is not of the directional smallest scale size.

According to Model (21), we know that we can find the optimal solutions of Model (21), denoted by $(\beta^*, \eta^*, \lambda_j^*)$, in which $\eta^* > 0$. Because $(X_0, Y_0)$ is strongly efficient, we know that $\beta^* > 0$. According to **Theorems 4 and 5**, we know that we can find a small enough positive constant $t_{left}$ that can ensure that the optimal value of Model (22) is constant. According to **Theorem 6**, we know that the upper bound $\bar{\rho}(X_0, Y_0)$ of Model (28) is equal to the optimal objective value $\psi^*(X_0, Y_0)$ of Model (22). This fact contradicts the supposition that the maximal optimal objective function $\bar{\rho}(X_0, Y_0)$ of Model (28) is unbounded. **Q.E.D.**